\documentclass[11pt]{amsart}
\usepackage{amsmath}
\usepackage{amsfonts}
\usepackage{latexsym}
\usepackage{amssymb}
\usepackage{enumerate}
\usepackage[dvips]{graphics}

\newcommand{\Z}{{\mathbf Z}}
\newcommand{\Q}{{\mathbf Q}}
\newcommand{\R}{{\mathbf R}}
\newcommand{\C}{{\mathbf C}}

\newcommand{\coker}{{\mbox{\rm coker}}}

\newcommand{\diag}{\rm {diag}}
\newcommand{\rk}{{\rm {rk}}}
\newcommand{\ind}{\rm {ind}}

\newtheorem{theorem}{Theorem}

\newtheorem{lemma}[theorem]{Lemma}
\newtheorem{corollary}[theorem]{Corollary}
\newtheorem{example}{Example}

\begin{document}
\title{Homology of planar polygon spaces}
\author[M. ~Farber and
D. ~Sch\"utz]{M.~Farber and D. ~Sch\"utz}
\address{Department of Mathematics, University of Durham, Durham DH1 3LE, UK}
\email{Michael.Farber@durham.ac.uk}

\email{Dirk.Schuetz@durham.ac.uk}




\date{}


\begin{abstract} In this paper we study topology of
the variety of closed planar $n$-gons with given side lengths
$l_1, \dots, l_n$. The moduli space $M_\ell$ where $\ell =(l_1,
\dots, l_n)$, encodes the shapes of all such $n$-gons. We describe
the Betti numbers of the moduli spaces $M_\ell$ as functions of
the length vector $\ell=(l_1, \dots, l_n)$. We also find sharp
upper bounds on the sum of Betti numbers of $M_\ell$ depending
only on the number of links $n$. Our method is based on an
observation of a remarkable interaction between Morse functions
and involutions under the condition that the fixed points of the
involution coincide with the critical points of the Morse
function.
\end{abstract}

\maketitle
\section{Introduction and statement of the result}

Given a string $\ell=(l_1, \dots, l_n)$ of $n$ positive real
numbers $l_i>0$ one considers the moduli space $M_\ell$ of closed
planar polygonal curves having side lengths $l_i$. Points of
$M_\ell$ parametrize different shapes of such polygons. Formally
$M_\ell$ is defined as the factor space
\begin{eqnarray*}
M_\ell = \{(u_1, \dots, u_n)\, \in \, S^1\times \dots\times S^1;
\, \sum_{i=1}^n l_iu_i \, =\, 0\in \C\}\, /\, {\rm {SO}}(2).
\end{eqnarray*}
Here $u_i\in S^1\subset \C$ denote the unit vectors in the
directions of the sides of a polygon; the group of rotations
$SO(2)$ acts diagonally on $(u_1, \dots, u_n)$.

Viewed differently, $M_\ell$ is the configuration space of a
planar linkage, a planar mechanism consisting of $n$ bars of
length $l_1, \dots, l_n$ connected by revolving joints. Such
mechanisms play an important role in robotics where they describe
closed kinematic chains and are used widely as elementary parts of
more complicated mechanisms. Knowing the topology of $M_\ell$ (for
different vectors $\ell$) can be used in designing control
programmes and motion planning algorithms for mechanisms.

The length vector $\ell$ is called {\it generic} if
$\sum\limits_{i=1}^n l_i\epsilon_i\not=0$ for any choice
$\epsilon_i=\pm 1$. It is known that for a generic length vector
$\ell$ the space $M_\ell$ is a closed smooth manifold of dimension
$n-3$. If the length vector $\ell$ is not generic then $M_\ell$ is
a compact $(n-3)$-dimensional manifold with finitely many singular
points.

The moduli spaces $M_\ell$ of planar polygonal linkages were
studied extensively by many mathematicians; we will mention W.
Thurston and J. Weeks \cite{TW}, K. Walker \cite{Wa},
 A. A. Klyachko \cite{Kl},
M. Kapovich and J. Millson \cite{KM1}, J.-Cl. Hausmann and A. Knutson \cite{HK}
and others.

Our goal in this paper is to give a general formula for the Betti
numbers of the moduli space $M_\ell$ as functions of the length
vector $\ell$. Our results cover both generic and non-generic
vectors $\ell$. In the case of generic $\ell$ the Betti numbers of
$M_\ell$ can easily be extracted from the results of the
unpublished thesis of K. Walker \cite{Wa}.

Formulae for Betti numbers of polygon spaces in three-dimensional
space are known (see papers of A.A. Klyachko \cite{Kl} and J.-Cl.
Hausmann and A. Knutson \cite{HK}). Also, J.-Cl. Hausmann and A.
Knutson  describe cohomology with $\Z_2$ coefficients of the
factor $\bar{M_\ell}$ of $M_\ell$ with respect to the natural
involution, see \cite{HK}, Theorem 9.1.

A. Klyachko in his beautiful work \cite{Kl} uses a remarkable
symplectic structure on the moduli space of linkages in $\R^3$ in
an essential way. His technique is based on properties of
Hamiltonian circle actions (the perfectness of the Hamiltonian
viewed as a Morse function). The subsequent important paper of
J.-Cl. Hausmann and A. Knutson employs methods of symplectic
topology as well: they apply the method of symplectic reduction.
Note that J.-Cl. Hausmann and A. Knutson go one step further and
compute the multiplicative structure on cohomology, however their
description is not very explicit as it uses the language of
generators and relations. Symplectic methods play also a central
role in the work of M. Kapovich and J. Millson \cite{KM2}.

The moduli spaces of planar linkages $M_\ell$ do not carry symplectic
structures in general. Therefore methods of symplectic topology are not
applicable in this problem.

The proof of our main result (see Theorem \ref{thm1} below) is obtained in a
very simple manner, it uses a remarkable interaction between Morse functions
and involutions under the condition that fixed points of the involution
coincide with the critical points of the Morse function.

To state our main theorem we need the following definitions. A
subset $J\subset \{1, \dots, n\}$ is called {\it short} if
$$\sum_{i\in J} l_i < \sum_{i\notin J} l_i.$$ The complement of a
short subset is called {\it long}. A subset $J\subset \{1, \dots,
n\}$ is called {\it median} if
$$\sum_{i\in J} l_i = \sum_{i\notin J} l_i.$$
Clearly, median subsets exist only if the length vector $\ell$ is
not generic. Note the following simple observation: any two
subsets $J, J'\subset \{1, \dots, n\}$ have a nonempty
intersection $J\cap J'\not=\emptyset$ provided that one of the
subsets is long and the other is either long or median.

\begin{theorem}\label{thm1} Fix a link of the maximal length $l_i$,
i.e. such that $l_i\geq l_j$ for any $j=1, 2, \dots, n$. For every
$k=0,1, \dots, n-3$ denote by $a_k$ and $b_k$ correspondingly the
number of short and median subsets of $\{1, \dots, n\}$ of
cardinality $k+1$ containing $i$. Then the homology group
$H_k(M_\ell;\Z)$ is free abelian of rank
\begin{eqnarray}\label{betti} a_k + b_k + a_{n-3-k},\end{eqnarray} for any $k=0, 1,
\dots, n-3$.
\end{theorem}

By Theorem \ref{thm1} the Poincar\'e polynomial $$p(t) = \sum_{k=0}^{n-3} \,
\dim H_k(M_\ell;\Q)\, \cdot\,  t^k$$ of $M_\ell$ can be written in the form
\begin{eqnarray} q(t) + t^{n-3}q(t^{-1}) + r(t)
\end{eqnarray} where
\begin{eqnarray}
q(t) =\sum_{k=0}^{n-3} a_k t^k, \quad r(t) = \sum_{k=0}^{n-3} b_k
t^k;
\end{eqnarray}
the numbers $a_k$ and $b_k$ are described in the statement of
Theorem \ref{thm1}.

A proof of Theorem \ref{thm1} is given below in \S \ref{proofthm1}. In the rest
of this introduction we illustrate the statement of Theorem \ref{thm1} by
several examples.

\begin{example} {\rm Suppose that $n=5$ and $l_1=3$, $l_2=2$, $l_3=2$, $l_4=1,
l_5=1$. Then $l_1=3$ is the longest link and short subsets of $\{1, \dots, 5\}$
containing $1$ are $\{1\}$, $\{1, 4\}$ and $\{1, 5\}$. Hence $a_0=1$, $a_1=2$
and by Theorem \ref{thm1} the Poinca\'re polynomial of $M_\ell$ equals
$1+4t+t^2$. We conclude that $M_\ell$ is a closed orientable surface of genus
2.}
\end{example}

\begin{example}\label{excon}{\rm Consider the zero-dimensional Betti number $$a_0+b_0+ a_{n-3}$$
of $M_\ell$ as given by Theorem \ref{thm1}. We want to show that
this number can take values $0, 1, 2$; the first possibility is
clearly equivalent to $M_\ell=\emptyset$. Without loss of
generality we may assume that $l_1\leq l_2\leq \dots\leq l_n$. If
$\{n\}$ is short then $a_0=1$ and $b_0=0$. If $\{n\}$ is median
then $a_0=0$ and $b_0=1$; in this case clearly $M_\ell$ is a
single point. If $\{n\}$ is long then $a_k=0=b_k$ for any $k$ and
hence $M_\ell=\emptyset$. We obtain that $M_\ell=\emptyset$ if and
only if there are no long one-element subsets of $\{1, \dots, n\}$
-- a result first established by Kapovich and Millson in
\cite{KM1}.

Let us show that the number $a_{n-3}$ equals $0$ or $1$. Clearly,
$a_{n-3}$ coincides with the number of long two-element subsets
$\{r, s\}\subset \{1, \dots, n-1\}$. There may exist at most one
such pair: if $\{r', s'\}$ is another long pair with $r\not= r'$,
$r\not=s'$, then $\{r, n\}$ and $\{r',s'\}$ would be two disjoint
long subsets which is impossible. We obtain that $a_{n-3}=1$ if
and only if the pair $\{n-2,n-1\}$ is long and $a_{n-3}=0$
otherwise.

We see that the moduli space $M_\ell$ has two connected components
if and only if the set $\{n-2, n-1\}$ is long and $\{n\}$ is
short. In this case the length vector $\ell$ must be generic and
short subsets $J\subset \{1, \dots, n\}$ containing $n$ are
exactly the subsets containing neither $n-2$ nor $n-1$. We see
that the Poincar\'e polynomial of $M_\ell$ in this case equals
$2(1+t)^{n-3}$. M. Kapovich and J. Millson \cite{KM1} showed that
if $M_\ell$ is disconnected then it is diffeomorphic to the
disjoint union of two copies of the torus $T^{n-3}$.}
\end{example}

\begin{example}\label{ex2}{\rm  As another example consider
the equilateral case when $l_j=1$ for all $j$. Assume first that
$n=2r+1$ is odd and hence $\ell$ is generic. The short subsets in
this case are subsets of $\{1, \dots, n\}$ of cardinality $\leq
r$. We may fix the index $\{n\}$ as representing the longest link.
Hence we find that $b_k=0$ vanishes and $a_k$ equals
\begin{eqnarray}\label{ak}
a_k= \left\{
\begin{array}{cl}
\left(\begin{array}{c} n-1\\ k\end{array}\right)& \mbox{for $k\leq r-1$,}\\
\\0,& \mbox{for $k\geq r$.}
\end{array}
\right.
\end{eqnarray}
By Theorem \ref{thm1} the Betti numbers of $M_\ell$ are given by
\begin{eqnarray}
b_k(M_\ell) = \left\{
\begin{array}{cl}
\left(\begin{array}{c} n-1\\ k\end{array}\right) & \mbox{for $k< r-1$},\\ \\

2\cdot \left(\begin{array}{c} n-1\\ r-1\end{array}\right) & \mbox{for
$k=r-1$},\\ \\

\left(\begin{array}{c} n-1\\ k+2\end{array}\right) & \mbox{for $k>r-1$}.
\end{array}
\right.
\end{eqnarray}
  Note that the sum of Betti
numbers in this example equals
\begin{eqnarray}\label{sb}
\sum_{k=0}^{n-3} b_k(M_\ell)\, = \, 2^{n-1}- \left(\begin{array}{c} n-1\\
r\end{array}\right), \quad \mbox{where}\quad n=2r+1.
\end{eqnarray}}\end{example}

\begin{example}\label{ex3}{\rm Consider now the equilateral case $l_j=1$ with
$n$ is even, $n=2r+2$. The length vector is now not generic. The
short subsets are all subsets of cardinality $\leq r$ and the
median subsets are all subsets of cardinality $r+1$. Hence we find
that $b_k=0$ for $k\not=r$ and
\begin{eqnarray}
b_r = \left(\begin{array}{c} 2r+1\\ r\end{array}\right)
\end{eqnarray}
and the numbers $a_k$ are given by formula (\ref{ak}). Applying
Theorem \ref{thm1} we find
\begin{eqnarray}
b_k(M_\ell) = \left\{
\begin{array}{cl}
\left(\begin{array}{c} n-1\\ k\end{array}\right) & \mbox{for $k\leq r-1$},\\ \\

 \left(\begin{array}{c} n\\ r\end{array}\right) & \mbox{for
$k=r$},\\ \\

\left(\begin{array}{c} n-1\\ k+2\end{array}\right) & \mbox{for
$r+1 \leq k\leq n-3$}.
\end{array}
\right.
\end{eqnarray}
The sum of Betti numbers in this example is
\begin{eqnarray}\label{sbeven1}
\sum_{k=0}^{n-3} b_k(M_\ell)\, = \, 2^{n-1}- \left(\begin{array}{c} n-1\\
r\end{array}\right), \quad \mbox{where}\quad n=2r+2.
\end{eqnarray}
}\end{example}

The results described in Examples \ref{ex2} and \ref{ex3} were
obtained earlier in \cite{Ka1}, \cite{Ka2} by different methods.

In the next section we shall see that Examples \ref{ex2} and
\ref{ex3} give moduli spaces $M_\ell$ with the maximal possible
total Betti number for all length vectors $\ell$ having the given
number of links $n$.

\section{Maximum of the total Betti number of $M_\ell$}

It is well known that the moduli space of pentagons $M_\ell$ with
a generic length vector $\ell=(l_1, \dots, l_5)$ is a compact
orientable surface of genus not exceeding $4$, see \cite{MT}. In
the equilateral case, i.e. if $\ell=(1,1,1,1,1)$, $M_\ell$ is
indeed an orientable surface of genus $4$ (it is a special case of
(\ref{sb})) and hence the above upper bound for pentagons is
sharp. In this section we state a theorem generalizing this result
for arbitrary $n$. Namely, we prove that for any length vector
$\ell=(l_1, \dots, l_n)$ the sum of the Betti numbers
\begin{eqnarray}\label{sum1}\sum_{i=0}^{n-3} b_i(M_{\ell})\end{eqnarray}
is less or equal than the sum of Betti numbers of the moduli space
of the equilateral linkage with the same number of sides $n$.

\begin{theorem}\label{upper} Let $\ell=(l_1, \dots, l_n)$ be a length vector, $l_i>0$.
 Denote by $r$ the number $[\frac{n-1}{2}]$.
 Then the sum of Betti numbers of the moduli space $M_\ell$ does not exceed
\begin{eqnarray}\label{sbm}
B_n = 2^{n-1} -\left(\begin{array}{c} n-1\\ r\end{array}\right).
\end{eqnarray}
This estimate is sharp: $B_n$ equals the sum of Betti numbers of
the moduli space of planar equilateral $n$-gons, see (\ref{sb}),
(\ref{sbeven1}).
\end{theorem}

Note that for $n$ even the equilateral linkage with $n$ sides is
not generic and hence Theorem \ref{upper} does not answer the
question about the maximum of the total Betti number on the set of
all generic length vectors with $n$ even.

\begin{theorem}\label{upper1} Assume that $n$ is even and $\ell=(l_1, \dots,
l_n)$ is a generic length vector. Then the sum of Betti numbers of
$M_\ell$ does not exceed
\begin{eqnarray}\label{sbmeven}
B'_n \, =\, 2\cdot B_{n-1},
\end{eqnarray}
where $B_{k}$ is defined by (\ref{sbm}). This upper bound is
achieved on the length vector $\ell=(1, 1, \dots, 1,\epsilon)$
where $0<\epsilon<1$ and the number of ones is $n-1$.
\end{theorem}

Note that $M_{(1,\dots, 1, \epsilon)}$ is diffeomorphic to the
product $M_{(1,\dots, 1)}\times S^1$ (the number of ones in both
cases equals $2r+1$). Hence the sum of Betti numbers of
$M_{(1,\dots, 1, \epsilon)}$ is twice the sum of Betti numbers of
$M_{(1,\dots, 1)}$.

Proofs of Theorems \ref{upper} and \ref{upper1} are given below in
section \S \ref{last}.

The asymptotic behavior of $B_n$ (given by (\ref{sbm})) can be
recovered using available information about Catalan numbers
$$C_r = \frac{1}{r+1} \cdot \left(
\begin{array}{c}
2r\\ r
\end{array}
\right) \sim \frac{2^{2r}}{\sqrt{\pi} r^{3/2}},
$$
see \cite{Vardi}. One obtains the following asymptotic formula
\begin{eqnarray}\label{bn}
B_n \sim 2^{n-1}\cdot\left( 1- \sqrt{\frac{2}{n\pi}}\right)
\end{eqnarray}
which is valid for even and odd $n$.

From the discussion of Example \ref{excon} we know that in the
case when $M_\ell$ is disconnected the sum of Betti numbers of
$M_\ell$ equals $2^{n-2}$ which is approximately half of $B_{n}$,
see (\ref{bn}).

\section{Morse theory on manifolds with involutions}

Our main tool in computing the Betti numbers of the moduli space of planar
polygons $M_\ell$ is Morse theory of manifolds with involution.

\begin{theorem}\label{thm2}
Let $M$ be a smooth compact manifold with boundary. Assume that $M$ is equipped
with a Morse function $f: M\to [0,1]$ and with a smooth involution $\tau: M\to
M$ satisfying the following properties:
\begin{enumerate}
\item $f$ is $\tau$-invariant, i.e. $f(\tau x)=f(x)$ for any $x\in M$; \item
The critical points of $f$ coincide with the fixed points of the involution;
\item $f^{-1}(1)=\partial M$ and $1\in [0,1]$ is a regular value of $f$.
\end{enumerate}
Then each homology group $H_i(M;\Z)$ is free abelian of rank equal the number
of critical points of $f$ having Morse index $i$. Moreover, the induced map
$$\tau_\ast: H_i(M;\Z) \to H_i(M;\Z)$$ coincides with multiplication by
$(-1)^i$ for any $i$.
\end{theorem}
\begin{figure}[h]
\resizebox{3.2cm}{4cm} {\includegraphics[161,
351][445,712]{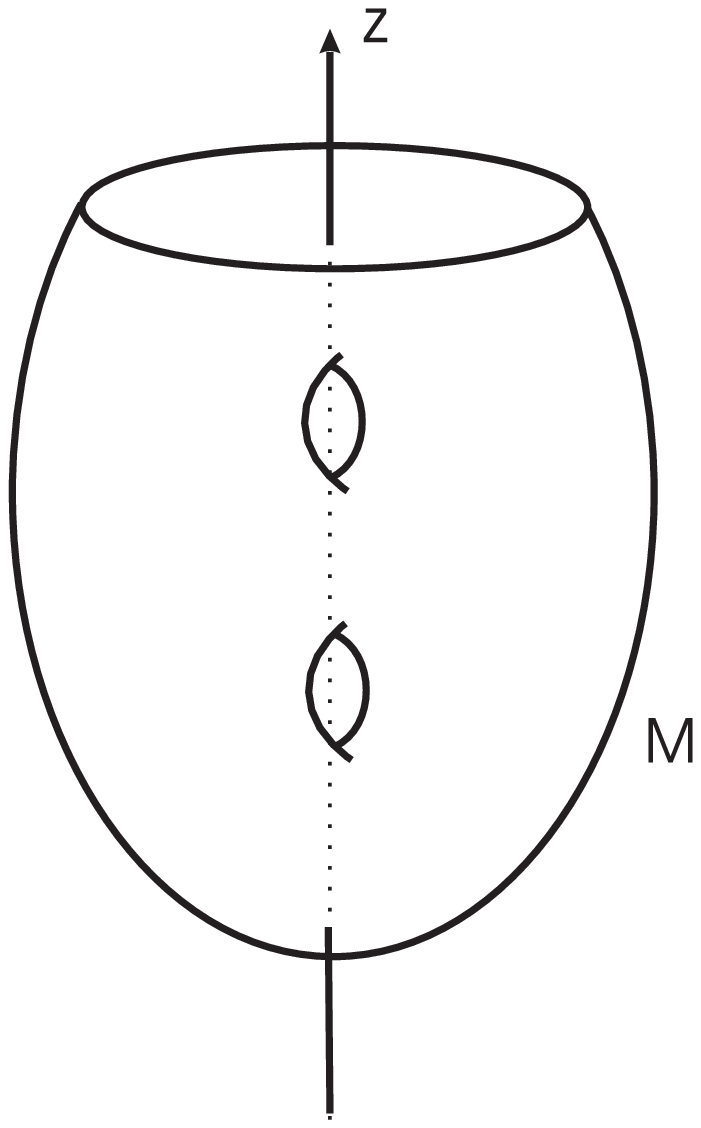}}
\caption{Surface in $\R^3$}\label{surface}
\end{figure}
As an illustration for Theorem \ref{thm2} consider a surface is $\R^3$ (see
Figure \ref{surface}) which is symmetric with respect to the
$z$-axis. The function $f$ is the orthogonal projection onto the
$z$-axis, the involution $\tau: M\to M$ is given by $\tau(x,
y,z)=(-x, -y, z)$.

The critical points of $f$ are exactly the intersection points of $M$ with the
$z$-axis.

\begin{proof}[Proof of Theorem \ref{thm2}]
Choose a Riemannian metric on $M$ which is invariant with respect to $\tau$.

Let $p\in M$ be a critical point of $f$. By our assumption, $p$ must be a fixed
point of $\tau$, i.e. $\tau(p)=p$. We claim that the differential of $\tau$ at
$p$ is multiplication by $-1$, i.e.
\begin{eqnarray}\label{minus}
d\tau_p(v) =-v, \quad \mbox{for any}\quad v\in T_pM.
\end{eqnarray}
Firstly, since $\tau$ is an involution, $d\tau_p$ must have eigenvalues $\pm
1$. Assume that there exists a vector $v\in T_pM$ with $d\tau_p(v)=v$. Then the
geodesic curve starting from $p$ in the direction of $v$ is invariant with
respect to $\tau$ implying that $p$ is not isolated in the fixed point set of
$\tau$. This contradicts our assumption and hence $d\tau_p$ must have
eigenvalue $-1$ only. Note that $d\tau_p$ is diagonalizable as $d\tau_p$
preserves the Hessian of $f$ at $p$
\begin{eqnarray} H(f)_p: T_p(M)\otimes T_p(M)\to \R
\end{eqnarray}
which is a nondegenerate quadratic form. This proves (\ref{minus}).

Consider the gradient vector field $v$ of $f$ with respect to the Riemannian
metric. We will assume that $v$ satisfies the transversality condition, i.e.
all stable and unstable manifolds of the critical points intersect
transversally. $v$ is $\tau$-invariant which means that
\begin{eqnarray}\label{inv}
v_{\tau(x)}=d\tau_x (v_x), \quad x\in M.
\end{eqnarray}

The Morse - Smale chain complex $(C_\ast(f),\partial)$ of $f$ has the critical
points of $f$ as its basis and the differential is given by
\begin{eqnarray}
\partial (p) = \sum_{q} [p:q]\, q
\end{eqnarray}
where in the sum $q$ runs over the critical points $q$ with Morse index ${\rm
{ind}}(q)= {\rm {ind}}(p)-1$. The incidence numbers $[p:q]\in \Z$ are defined
as follows
\begin{eqnarray}\label{incid}
[p:q] = \sum_\gamma \epsilon(\gamma), \quad 
\epsilon(\gamma)=\pm 1,
\end{eqnarray}
where $\gamma: (-\infty, \infty)\to M$ are trajectories of the negative
gradient flow $\gamma'(t) = -v_{\gamma(t)}$ satisfying the boundary conditions
$\gamma(t)\to p$ as $t\to -\infty$ and $\gamma(t)\to q$ as $t\to +\infty$.

Observe that if $\gamma$ is a trajectory as above then $\tau\circ \gamma$ is
another such trajectory. Indeed, using (\ref{inv}) we find $(\tau\circ \gamma)'
=d\tau(\gamma') = -d\tau (v_{\gamma(t)}) = -v_{\tau(\gamma(t))}.$

 Theorem \ref{thm2} would follow once we show that
\begin{eqnarray}\label{zero}
\epsilon (\gamma) + \epsilon(\tau\circ \gamma) =0,\end{eqnarray}
i.e. the total contribution into (\ref{incid}) of a pair of
symmetric trajectories is zero. Hence all incidence coefficients
vanish $[p:q]=0$ and the differentials of the Morse - Smale
complex are trivial.\begin{figure}[h]
\resizebox{4cm}{3.8cm}
{\includegraphics[164,385][512,710]{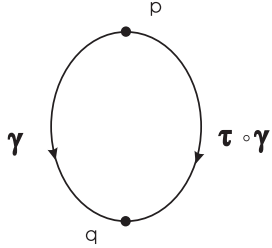}}
\caption{Two symmetric trajectories of the negative gradient
flow}\label{two}\end{figure}

To prove (\ref{zero}) we first recall the definition of the sign
$\epsilon(\gamma)\in \{1, -1\}$, see \cite{Mi2}. For a critical point $p$ of
$f$ we denote by $W^u(p)$ and $W^s(p)$ the unstable and stable manifolds of
$p$. Recall that $W^u(p)$ is the union of the trajectories $\gamma: (-\infty,
\infty)\to M$ satisfying the differential equation
 $\gamma'(t) = -v_{\gamma(t)}$ and the boundary condition $\gamma(t)\to p$ as
$t\to -\infty$. The stable manifold $W^s(p)$ is defined similarly but the
boundary condition in this case becomes $\gamma(t)\to p$ as $t\to +\infty$.

Fix an orientation of the stable manifold $W^s(p)$ for every critical point
$p\in M$. Since $W^s(p)$ and $W^u(p)$ are of complementary dimension and
intersect transversally at $p$, the orientation of $W^s(p)$ determines a
coorientation of the unstable manifold $W^u(p)$, for every $p$.

If $\ind(p) -\ind(q) =1$
 then $W^u(p)$ and $W^s(q)$ intersect
transversally along finitely many connecting orbits $\gamma(t)$ and the
structure near each of the connecting orbits looks as shown on Figure
\ref{intersect}.
\begin{figure}[h]
\resizebox{5cm}{5cm} {\includegraphics[94,304][464, 739]{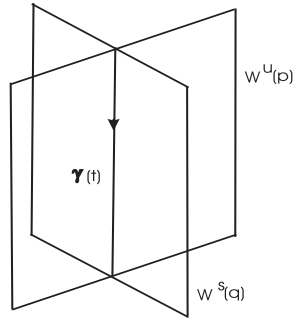}}
\caption{The stable and unstable manifolds along $\gamma(t)$}\label{intersect}
\end{figure}
Note that the normal bundle to $W^u(p)$ along $\gamma$ coincides with the
normal bundle to $\gamma$ in $W^s(q)$. Hence, the coorientation of $W^u(p)$
together with the natural orientation of the curve $\gamma(t)$ determine an
orientation of $W^s(q)$ along $\gamma$. We set $\epsilon(\gamma)= 1$ iff this
orientation coincides with the prescribed orientation of $W^s(q)$; otherwise we
set $\epsilon(\gamma)=-1$.

To compare $\epsilon(\gamma)$ with $\epsilon(\tau\circ \gamma)$ we first
observe that the involution $\tau$ preserves the stable and unstable manifolds
$W^s(p)$ and $W^u(p)$ and for every critical point $p$ the degrees of the
restriction of $\tau$ on these submanifolds equal
\begin{eqnarray}\label{minus1}
\deg(\tau|_{W^u(p)})=(-1)^{\ind(p)}, \quad \quad
\deg(\tau|_{W^s(p)})=(-1)^{n-\ind(p)},
\end{eqnarray}
as follows from (\ref{minus}). Hence, applying the involution $\tau$ to the
picture shown on Figure \ref{intersect}, we have to multiply the coorientation
of $W^u(p)$ by $(-1)^{n-i-1}$ and multiply the orientation of $W^s(q)$ by
$(-1)^{n-i}$. As the result the total sign will be multiplied by
$(-1)^{n-i-1}\cdot (-1)^{n-i}= -1$. This proves (\ref{zero}) and completes the
proof of the first statement of the theorem. The second statement of the
Theorem follows from the first one combined with (\ref{minus1}).
\end{proof}

\begin{theorem}\label{thm3}
Let $M$ be a smooth compact connected manifold with boundary. Suppose that $M$
is equipped with a Morse function $f: M\to [0,1]$ and with a smooth involution
$\tau: M\to M$ satisfying the properties of Theorem \ref{thm2}. Assume that for
any critical point $p\in M$ of the function $f$ we are given a smooth closed
connected submanifold $$X_p\subset M$$ with the following properties:
\begin{enumerate}
\item $X_p$ is $\tau$-invariant, i.e. $\tau(X_p)=X_p$; \item $p\in X_p$ and for
any $x\in X_p-\{p\}$, one has $f(x) < f(p)$; \item the function $f|_{X_p}$ is
Morse and the critical points of the restriction $f|_{X_p}$ coincide with the
fixed points of $\tau$ lying in $X_p$. In particular, $\dim X_p=\ind (p)$.
\item For any fixed point $q\in X_p$ of $\tau$ the Morse indexes of $f$ and of
$f|_{X_p}$ at $q$ coincide.
\end{enumerate}
Then each submanifold $X_p$ is orientable and the set of homology classes
realized by $\{X_p\}_{p\in {\rm {Fix}}(\tau)}$ forms a free basis of the
integral homology group $H_\ast(M;\Z)$. In other words, we claim that the
inclusion induces an isomorphism
\begin{eqnarray}
\bigoplus_{\ind(p)=i} H_i(X_p;\Z) \, \to\,  H_i(M;\Z)
\end{eqnarray}
for any $i$.
\end{theorem}
\begin{proof}[Proof of Theorem \ref{thm3}] First we note that each submanifold
$X_p$ is orientable. Indeed, Theorem \ref{thm2} applied to the restriction
$f|_{X_p}$ implies that $f|_{X_p}$ has a unique maximum and unique minimum and
the top homology group $H_i(X_p;\Z)=\Z$ is infinite cyclic where $i=\dim
X_p=\ind(p)$.

For a regular value $a\in \R$ of $f$ we denote by $M^a\subset M$ the preimage
$f^{-1}(-\infty, a]$. It is a compact manifold with boundary. It follows from
Theorem \ref{thm2} that $f$ has a unique local minimum and therefore $M^a$ is
either empty or connected. For $a$ slightly above the minimum value $f(p_0) =
\min f(M)$ the manifold $M^a$ is a disc and the homology of $M^a$ is obviously
realized by the submanifold $X_{p_0}\subset M^a$.

We proceed by induction on $a$. Our inductive statement is that the homology of
$M^a$ is freely generated by the homology classes of the submanifolds $X_p$
where $p$ runs over all critical points of $f$ satisfying $f(p)\leq a$.

Suppose that the statement is true for $a$ and the interval $[a, b]$ contains a
single critical value $c$. Let $p_1, \dots, p_r$ be the critical points of $f$
lying in $f^{-1}(c)$. Denote
$$X = \coprod_{i=1}^r X_{p_i}$$
(the disjoint union). Then $f$ induces a Morse function $\bar f: X\to \R$ and
we set $$X^a= \bar f^{-1}(-\infty, a].$$ Consider the Morse - Smale complexes
$C_\ast(M^a)$, $C_\ast(M^{b})$, $C_\ast(X)$ and $C_\ast(X^a)$; the first two
are constructed using the function $f$ and the latter two are constructed using
the function $\bar f$. We have the following Mayer-Vietoris-type short exact
sequence of chain complexes
\begin{eqnarray} \label{morsechain} 0\to C_\ast(X^a) \to C_\ast(X) \oplus
C_\ast(M^a) \stackrel \Phi\to C_\ast(M^{b})\to 0
\end{eqnarray} which (by the arguments indicated in the proof of Theorem
\ref{thm2}) have trivial differentials and hence the sequence
\begin{eqnarray}\label{MV}
0\to H_i(X^a) \to H_i(X)\oplus H_i(M^a) \stackrel \Phi\to H_i(M^{b})\to 0
\end{eqnarray}
is exact (all homology groups have coefficients $\Z$). It follows from Lemma
\ref{4} below and the construction of the Morse - Smale complex (compare
\cite{Mi2}, \S 7) that the homomorphism $\Phi$ (which appears in
(\ref{morsechain}) and (\ref{MV})) coincides with the sum of the chain maps
induced by the inclusions $X\to M^b$ and $M^a\to M^b$.

For $i< \dim X$ we have $H_i(X^a) \to H_i(X)$ is an isomorphism (by Theorem
\ref{thm2}) and hence (\ref{MV}) implies that $H_i(M^a)\to H_i(M^{b})$ is an
isomorphism. For $i\geq \dim X$ we have $H_i(X^a)=0$ and therefore $\Phi:
H_i(X)\oplus H_i(M^a)\to H_i(M^{b})$ is an isomorphism. This completes the step
of induction.
\end{proof}

Here is a minor variation of the Morse lemma which has been used in the proof.

\begin{lemma}\label{4} Let $f: \R^n\to \R$ be a smooth function having $0\in \R^n$ as a
nondegenerate critical point and suppose that for some $k\leq n$ the
restriction $f|_{\R^k\times \{0\}}: \R^k\times \{0\} \to \R$ also has a
nondegenerate critical point at $0\in \R^k$. Then there exists a neighborhood
$U\subset \R^n$ of $0$ and a local coordinate system $x: U\to \R^n$ such that
$x(\R^k\times\{0\}\cap U) \subset \R^k\times \{0\}$ and
\begin{eqnarray}
f(x_1, \dots, x_n) = \pm  x_1^2 + \dots + \pm x_n^2 + f(0).
\end{eqnarray}
\end{lemma}
\begin{proof} One simply checks that the coordinate changes in the standard
proof of the Morse lemma (compare \cite{Mi1}, \S 2) can be chosen so that the
subspace $\R^k\times \{0\}$ is mapped to itself.
\end{proof}

\section{The robot arm distance map}

\begin{figure}[h]
\resizebox{7cm}{5cm} {\includegraphics[90,368][440,632]{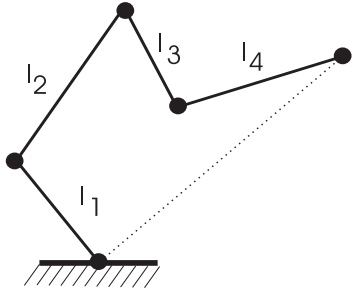}}
\caption{Robot arm.}\label{arm}
\end{figure}

A robot arm is a simple mechanism consisting of $n$ bars (links) of fixed
length $(l_1, \dots, l_n)$ connected by revolving joints, see Figure \ref{arm}.
The initial point of the robot arm is fixed on the plane.
The moduli space of a robot arm (i.e. the space of its possible shapes) is
\begin{eqnarray}\label{w}
W=\{(u_1, \dots, u_n)\in S^1\times \dots\times S^1\}/{\rm {SO}}(2).
\end{eqnarray}
Clearly, $W$ is diffeomorphic to a torus $T^{n-1}$ of dimension
$n-1$. A diffeomorphism can be specified, for example, by assigning to a
configuration $(u_1, \dots, u_n)$ the point $(1, u_2u_1^{-1}, u_3u_1^{-1},
\dots, u_{n-1}u_1^{-1})\in T^{n-1}$ (measuring angles between the directions of
the first and the other links).

Consider the moduli space of polygons $M_\ell$ (where $\ell=(l_1, \dots, l_n)$)
which is naturally embedded into $W$.

We define a function on $W$ as follows:
\begin{eqnarray}\label{fl}
f_\ell: W\to \R,\quad \quad f_\ell(u_1, \dots, u_n)= - \left|\sum_{i=1}^n l_i
u_i\right|^2.
\end{eqnarray}
Geometrically the value of $f_\ell$ equals the negative of the squared distance
between the initial point of the robot arm to the end of the arm shown by the
dotted line on Figure \ref{arm}. Note that the maximum of $f_\ell$ is achieved
on the moduli space of planar linkages $M_\ell\subset W$.

\begin{figure}[h]
\resizebox{5cm}{2cm} {\includegraphics[87,477][475,643]{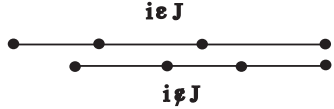}}
\caption{A collinear configuration $p_J$ of the robot arm.}\label{long}
\end{figure}

An important role play the collinear configurations, i.e. such
that $u_i=\pm u_j$ for all $i, j$, see Figure \ref{long}. We will
label such configurations by long and median subsets $J\subset
\{1, \dots, n\}$ assigning to any such subset $J$ the
configuration $p_J\in W$ given by $p_J=(u_1, \dots, u_n)$ where
$u_i=1$ for $i\in J$ and $u_i=-1$ for $i\notin J$. Note that $p_J$
lies in $M_\ell\subset W$ if and only if the subset $J$ is median.

\begin{lemma}\label{lemma1} The critical points of $f_\ell: W\to \R$ lying in
$W-M_\ell$ are exactly the collinear configurations $p_J$
corresponding to long subsets $J\subset \{1, 2, \dots, n\}$. Each
$p_J$, viewed as a critical point of $f_\ell$, is nondegenerate in
the sense of Morse and its Morse index equals $n-|J|$.
\end{lemma}

This lemma is well-known. It can be found as Proposition 3.3 in
\cite{Wa} and as combination of Theorems 3.1 and 3.2 in
\cite{Hau}; in both these references slightly different notations
were used.

\section{Proof of Theorem \ref{thm1}.}\label{proofthm1}

Consider the moduli space $W$ of the robot arm (defined by
(\ref{w})) with the function $f_\ell: W\to \R$ (defined by
(\ref{fl})). There is an involution
\begin{eqnarray}
\tau: W\to W
\end{eqnarray}
given by
\begin{eqnarray}\label{tau} \tau(u_1, \dots, u_n) = (\bar u_1, \dots, \bar u_n).
\end{eqnarray}
Here the bar denotes complex conjugation, i.e. the reflection with
respect to the real axis. It is obvious that formula (\ref{tau})
maps ${\rm {SO}}(2)$-orbits into ${\rm {SO}}(2)$-orbits and hence
defines an involution on $W$. The fixed points of $\tau$ are the
collinear configurations of the robot arm, i.e. the critical
points of $f_\ell$ in $W-M_\ell$, see Lemma \ref{lemma1}. Our plan
it to apply Theorems \ref{thm2} and \ref{thm3} to the sublevel
sets
\begin{eqnarray}\label{wa}
W^a=f_\ell^{-1}(-\infty,a]
\end{eqnarray}
of $f_\ell$. Recall that the values of $f_\ell$ are nonpositive
and the maximum is achieved on the submanifold $M_\ell\subset W$.
From Lemma \ref{lemma1} we know that the critical points of
$f_\ell$ are the collinear configurations $p_J$. The latter are
labelled by long subsets $J\subset \{1, \dots, n\}$ and $p_J=(u_1,
\dots, u_n)$ where $u_i=1$ for $i\in J$ and $u_i=-1$ for $i\notin
J$. One has \begin{eqnarray} f_\ell(p_J) = -(L_J)^2.\end{eqnarray}
Here $L_J=\sum_{i=1}^n l_iu_i$ with $p_J=(u_1, \dots, u_n)$.

\begin{figure}[h]
\resizebox{5cm}{4cm} {\includegraphics[51,304][652,750]{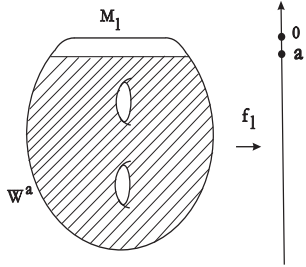}}
\caption{Function $f_\ell: W\to \R$ and the manifold
$W^a$.}\label{wa1}
\end{figure}

The number $a$ which appears in (\ref{wa}) will be chosen so that
\begin{eqnarray}\label{aa}
-(L_J)^2 < a < 0
\end{eqnarray}
for any long subset $J$ such that the manifold $W^a$ contains all
the critical points $p_J$. The situation is shown schematically on
Figure \ref{wa1}.

For each subset $J\subset \{1, \dots, n\}$ we denote by $\ell_J$
the length vector obtained from $\ell=(l_1, \dots, l_n)$ by
integrating all links $l_i$ with $i\in J$ into one link. For
example, if $J=\{1, 2\}$ then $\ell_J= (l_1+l_2, l_3, \dots,
l_n)$. We denote by $W_J$ the moduli space of the robot arm with
the length vector $\ell_J$. It is obvious that $W_J$ is
diffeomorphic to a torus $T^{n-|J|}$. We view $W_J$ as being
naturally embedded into $W$. Note that the submanifold $
W_J\subset W $ is disjoint from $M_\ell$ (in other words, $W_J$
contains no closed configurations) if and only if the subset
$J\subset \{1, \dots, n\}$ is long.

\begin{lemma}\label{lemma2}
Let $J\subset \{1, \dots, \}$ be a long subset. The submanifold
$W_J\subset W$ has the following properties:
\begin{enumerate}
\item $W_J$ is invariant with respect to the involution $\tau:
W\to W$; \item the restriction of $f_\ell$ onto $W_J$ is a Morse
function having as its critical points the collinear
configurations $p_I$ where $I$ runs over all subsets $I\subset
\{1, \dots, n\}$ containing $J$. \item for any such $p_I$ the
Morse indexes of $f_\ell$ and of $f_\ell|_{W_J}$ at $p_I$
coincide.

 \item in particular, $f|_{W_J}$ achieves its maximum at $p_J\in W_J$.
 \end{enumerate}
\end{lemma}
\begin{proof} (1) is obvious. Statements (2) and (3) follow from Lemma
\ref{lemma1} applied to the restriction of $f_\ell$ onto $W_J$.
Here we use the assumption that $J$ is long. Under this assumption
the long subset for the integrated length vector $\ell_J$ are in
one-to-one correspondence with the long subsets $I\subset \{1,
\dots,n\}$ containing $J$. Statement (4) follows from (3) as the
Morse index of $f_\ell|_{W_J}$ at point $p_J$ equals $n-|J|=\dim
W_J$.
\end{proof}

Applying Theorems \ref{thm2} and \ref{thm3} and taking into
account Lemma \ref{lemma2} we obtain:

\begin{corollary} One has:
\begin{enumerate}
\item If $a$ satisfies (\ref{aa}) then the manifold $W^a$ (see
(\ref{wa})) contains all submanifolds $W_J$ where $J\subset \{1,
\dots, n\}$ is an arbitrary long subset. \item The homology
classes of the submanifolds $W_J$ form a free basis of the
integral homology group $H_\ast(W^a;\Z)$.
\end{enumerate}
\end{corollary}

Next we examine the homomorphism
\begin{eqnarray}\label{phi}
\phi_\ast: H_i(W^a;\Z) \to H_i(W;\Z)
\end{eqnarray}
induced by the inclusion $\phi: W^a\to W$.

Below we will assume that $l_1\geq l_j$ for all $j\in \{1, \dots,
n\}$, i.e. $l_1$ is the longest link. This may always be achieved
by relabelling.

We describe a specific basis of the homology $H_\ast(W;\Z)$. For
any subset $J\subset \{1, 2, \dots, n\}$ we denote by $W_J$ the
moduli space of configurations of the robot arm with length vector
$\ell_J$ where all links $l_i$ with $i\in J$ are integrated into a
single link. Note that $W_J$ is naturally embedded into $W$ and
$$W_J\cap M_\ell=\emptyset$$ if and only if the set $J$ is long.
Since $W$ is homeomorphic to the torus $T^{n-1}$, it is easy to
see that a basis of the homology group $H_\ast(W;\Z)$ is formed by
the homology classes of the submanifolds $W_J$ where $J\subset
\{1, \dots, n\}$ runs over all subsets containing $1$. We will
denote the homology class of $W_J$ by \begin{eqnarray}\label{wi}
[W_J]\in H_{n-|J|}(W;\Z).
\end{eqnarray}

Assuming that $J, J'\subset \{1, \dots, n\}$ are two subsets with
$|J|+|J'| =n+1$ the classes $[W_J]$ and $[W_{J'}]$ have
complementary dimensions in $W$ and their intersection number is
given by
\begin{eqnarray}\label{intersection}
[W_J]\cdot [W_{J'}] = \left\{
\begin{array}{lll}
\pm 1, &\mbox{if} & |J\cap J'|=1, \\ \\
0, &\mbox{if} & |J\cap J'|>1.
\end{array} \right.
\end{eqnarray}
Indeed, if $J\cap J'=\{i_0\}$ then $W_J\cap W_{J'}$ consists of a
single point $\{p\}$, the moduli space of a robot arm with all
links integrated into one link. Let us show that the intersection
$W_J\cap W_{J'}$ is transversal. A tangent vector to $W$ at
$p=(u_1, \dots, u_n)$ can be labelled by a vector $w=(\lambda_1,
\dots, \lambda_n)\in\R^n$ (an element of the Lie algebra of the
torus $T^n$) viewed up to adding vectors of the form $(\lambda,
\lambda, \dots, \lambda)$. Such a tangent vector $w$ is tangent to
the submanifold $W_J$ iff $\lambda_i=\lambda_j$ for all $i, j\in
J$. Given $w$ as above it can be written as
$$w=w'+w''+(\lambda_{i_0}, \dots, \lambda_{i_0})$$
where $w'$ has coordinates $0$ on places $i\in J$ and coordinates
$\lambda_i-\lambda_{i_0}$ on places $i\not\in J$; coordinates of
$w''$ vanish on places $i\notin J$ and are
$\lambda_i-\lambda_{i_0}$ on places $i\in J$. Hence every tangent
vector to $W$ is a sum of a tangent vector to $W_J$ and a tangent
vector to $W_{J'}$.

Now suppose that $|J\cap J'|>1$. We will show that then the
submanifold $W_{J'}$ can be continuously deformed inside $W$ to a
submanifold $W'_{J'}$ such that $W_J\cap W'_{J'}=\emptyset$. This
would prove the second claim in (\ref{intersection}). Let us
assume that $\{1, 2\}\subset J\cap J'$. Define $g_t: W_{J'}\to W$
by $$g_t(u_1, \dots, u_n) = (e^{i\theta t}u_1, u_2, \dots, u_n),
\quad t\in [0,1].$$ Here $\theta$ satisfies $0<\theta<\pi$. Then
$W'_{J'}=g_1(W_{J'})$ is clearly disjoint from $W_J$; indeed, the
links $l_1$ and $l_2$ are parallel in $W_J$ and make an angle
$\theta$ in $W'_{J'}$.

It follows that the intersection form in the basis $[W_J],
[W_{J'}]\in H_\ast(W;\Z)$, where $J\ni 1, J'\ni 1$, has a very
simple form:
\begin{eqnarray}\label{intersection1}
[W_J]\cdot [W_{J'}] = \left\{
\begin{array}{lll}
\pm 1, &\mbox{if} & J\cap J'=\{1\}, \\ \\
0, &\mbox{if} & |J\cap J'|>1.
\end{array} \right.
\end{eqnarray}
In particular, given $[W_J]$ with $1\in J$, its dual homology
class $\in H_\ast(W;\Z)$ equals $[W_K]$ where $K=CJ\cup \{1\}$;
here $CJ$ denotes the complement of $J$ in $\{1, \dots, n\}$.

Denote by $A_\ast\subset H_\ast(W^a;\Z)$ (correspondingly,
$B_\ast\subset H_\ast(W^a;\Z)$) the subgroup generated by the
homology classes $[W_J]$ where $J\subset \{1, \dots, n\}$ is long
and contains $1$ (correspondingly, $J$ is long and $1\notin J$).
Then
\begin{eqnarray}\label{phi1}
H_i(W^a;\Z)= A_i\oplus B_i. \end{eqnarray} Similarly, one has
\begin{eqnarray}\label{phi2}
H_i(W;\Z)= A_i\oplus C_i \oplus D_i, \end{eqnarray} where:
\begin{itemize}
\item $A_\ast$ is as above; \item $C_\ast\subset H_\ast(W;\Z)$ is
the subgroup generated by the homology classes $[W_J]$ with
$J\subset \{1, \dots, n\}$ short and $1\in J$; \item $D_\ast$ is
the subgroup generated by the classes $[W_J]\in H_\ast(W; \Z)$
where $J$ is median and contains $1$.\end{itemize}

It is clear that $\phi_\ast$ (see (\ref{phi})) is identical when
restricted to $A_i$, compare (\ref{phi1}) and ({\ref{phi2}). We
claim that the image $\phi_\ast(B_i)$ is contained in $A_i$. This
would follow once we show that
\begin{eqnarray}\label{zero1}
[W_J]\cdot [W_{K}]=0
\end{eqnarray}
assuming that $[W_J]\in B_i$ and $[W_K]$ is the dual of a class
$[W_{J'}]\in C_i$ or $[W_{J'}]\in D_i$, see (\ref{intersection1}).
We have
\begin{enumerate} \item $J$ is long and $1\notin J$, \item $J'$
is short or median and $1\in J'$, \item $|J|= |J'|$, \item
$K=CJ'\cup \{1\}$.
\end{enumerate}
Here $CJ'$ denotes the complement of $J'$ in $\{1, \dots, n\}$. By
(\ref{intersection}), to prove (\ref{zero1}) we have to show that
under the above conditions one has $|J\cap K|>1$. Indeed, suppose
that $|J\cap K|=1$, i.e. $J\cap K=\{j\}$, a single element subset.
Then $J'$ is obtained from $J$ by removing the index $j$ and
adding the index $1$ which leads to a contradiction: indeed, $J$
is long, $l_j\leq l_1$ and $J'$ is either short or median.

\begin{corollary}
The kernel of the homomorphism $$\phi_i: H_i(W^a;\Z)\to
H_i(W;\Z)$$ has rank equal\footnote{Note that the kernel of
$\phi_i$ (viewed as a subgroup) is distinct from $B_i$ in
general.} to $\rk\, {B_i}$ and the cokernel has rank $\rk\, {C_i}+
\rk \, D_i$.
\end{corollary}

Below we skip the coefficient group $\Z$ from the notations.

One has
\begin{eqnarray}
\qquad H_{j}(W, W^a) \simeq H_{j}(N, \partial N)\simeq
H^{n-1-j}(N) \simeq H^{n-1-j}(M_\ell).
\end{eqnarray}
Here $N$ denotes the preimage $f_\ell^{-1}([a,0])$.

Note that $M_\ell$ is a deformation retract of $N$. Indeed,
consider $M_\ell\subset N^{\prime\prime}\subset N'\subset N$ where
$N'$ is a regular neighborhood of $M_\ell$ in $N$ and
$N^{\prime\prime}$ is a sublevel set
$N^{\prime\prime}=f_\ell^{-1}([a',0])$ and $a<a'<0$ is such that
$N^{\prime\prime}$ is contained in $N'$. Since $M_\ell\subset N'$
and $N^{\prime\prime}\subset N$ are deformation retracts, we have
the following diagram
$$
\begin{array}{ccc}
M_\ell & \stackrel {r'}\leftarrow & N'\\ \\
i\downarrow & \nearrow j & \downarrow k\\ \\
N^{\prime\prime} &
\stackrel r\leftarrow  & N
\end{array}
$$
where $i, j, k$ are inclusions and $r'ji=1_{M_\ell}$, $jir'\simeq
1_{N'}$, $rkj = 1_{N^{\prime\prime}}$, $kjr \simeq 1_N.$ It
follows  that $g=r'jr: N\to M_\ell$ is a deformation retraction.

Hence we obtain the following short exact sequence
\begin{eqnarray}\label{ssec}
0\to \coker (\phi_{n-1-j}) \to H^{j}(M_\ell) \to \ker
\phi_{n-2-j}\to 0
\end{eqnarray}
which splits since the kernel of $\phi_{n-2-j}$ is isomorphic to
$B_{n-2-j}$ (see above) and hence it is free abelian.

This proves that the cohomology $H^\ast(M_\ell)$ has no torsion
and therefore the homology $H_\ast(M_\ell)$ is free as well (by
the Universal Coefficient Theorem). The cokernel of $\phi_{n-1-j}$
is isomorphic to $C_{n-1-j}\oplus D_{n-1-j}$ as we established
earlier. We find that the rank of $\coker \phi_{n-1-j}$ equals the
number of subsets $J\subset\{1, \dots, n\}$ which are short or
median and have cardinality $|J|=j+1$. In other words,
\begin{eqnarray}\label{rkcoker}
\rk (\coker \phi_{n-1-j}) = a_j +b_j,
\end{eqnarray}
where we use the notation introduced in the statement of Theorem
\ref{thm1}.

The rank of the kernel of $\phi_{n-2-j}$ equals the rank of
$B_{n-2-j}$, i.e. the number of long subsets $J\subset \{2, \dots,
n\}$ of cardinality $|J|=j+2$. Passing to the complements, we find
\begin{eqnarray}\label{rkker} \rk (\ker \phi_{n-2-j})= a_{n-3-j} \end{eqnarray}
i.e. the number of short subsets containing $1$ with $|J|= n-2-j$.

Combining (\ref{rkcoker}), (\ref{rkker}) with the exact sequence
(\ref{ssec}) we finally obtain
$$\rk \, H_j(M_\ell) =\rk\,  H^j(M_\ell) = a_j +b_j + a_{n-3-j}.$$
This completes the proof, compare (\ref{betti}).

\section{Proofs of Theorems \ref{upper} and \ref{upper1}}\label{last}

The proofs are based on Theorem \ref{thm1} and are purely
combinatorial. Let $\ell=(l_1, \dots, l_n)$ be a length vector.
Without loss of generality we may assume that $l_1\leq l_2\leq
\dots\leq l_n$. By Theorem \ref{thm1} the sum of Betti numbers of
$M_\ell$ equals twice the number of short subsets plus the number
of median subsets of $\{1, \dots, n\}$, containing $n$. We show
that the number of such subsets is bounded above by $B_n$ (given
by (\ref{sbm})); moreover, we show that it is bounded above by
$B'_n=2\cdot B_{n-1}$ (see (\ref{sbmeven})) assuming additionally
that $\ell$ is generic and $n$ is even.

We will treat simultaneously both cases $n$ even and $n$ odd.
Denote $r = [(n-1)/2]$ so that $n=2r+2$ for $n$ even and $n=2r+1$
for $n$ odd.

For $1\leq i\leq n$ we denote by $S_i(\ell)$ (respectively,
$M_i(\ell)$) the number of short (respectively, median) subsets of
$\{1, \dots, n\}$ containing the subset $\{n-i+1, n-i+2, \dots,
n\}$. Clearly, $S_{r+1}(\ell) =M_{r+1}(\ell) =0$ for $n$ odd and
$S_{r+1}(\ell) = 0$, $M_{r+1}(\ell) \leq 1$ for $n$ even.

We claim that
\begin{eqnarray}\label{intermid}
\qquad 2\cdot S_i(\ell) + M_i(\ell)\, \leq\,  2^{n-i}\, - \,
\sum_{j=r-i+1}^r \left(\begin{array}{c} n-i\\
j\end{array}\right)
\end{eqnarray}
for all $1\leq i\leq r+1$. For $i=1$ inequality (\ref{intermid})
gives
$$2\cdot S_1(\ell) + M_1(\ell) \, \leq\,  2^{n-1} - \left(\begin{array}{c} n-1\\
r\end{array}\right) \, =\,  B_n$$ which is equivalent to our goal
(\ref{sbm}). We will prove (\ref{intermid}) by induction on $n$
and by descending induction on $i$.

For $n$ odd and $i=r+1$ inequality (\ref{intermid}) gives $2\cdot
S_{r+1}(\ell) + M_{r+1}(\ell) \leq 0,$ which follows from our
remark above. Similarly, for $n$ even and $i=r+1$ inequality
(\ref{intermid}) states
 $2\cdot S_{r+1}(\ell) + M_{r+1}(\ell) \leq 1,$ which is
obviously true, see above. These two remarks serve as the initial
step of induction.

Assume now that inequality (\ref{intermid}) is true (a) for $i+1$
and (b) for all $i$ and all length vectors $\ell'=(l'_1, \dots,
l'_m)$ with $m<n$ odd.

One can write
\begin{eqnarray}\label{adding}  S_i(\ell) =
S_{i+1}(\ell) + S'_i(\ell), \quad M_i(\ell) = M_{i+1}(\ell) +
M'_i(\ell)\end{eqnarray} where $S'_i(\ell)$ and $M'_i(\ell)$
denote the numbers of short and median subsets of $\{1, \dots,
n\}$ containing $\{n-i+1, \dots, n\}$ and not containing $n-i$.
One observes that
\begin{eqnarray}\label{fine}
\quad S'_i(\ell) \leq S_{i-1}(\tilde \ell), \quad \mbox{and}\quad
S'_i(\ell) + M'_i(\ell) \, \leq\, S_{i-1}(\tilde \ell)+
M_{i-1}(\tilde \ell),
\end{eqnarray}
where  \begin{eqnarray} \tilde \ell=(l_1, l_2, \dots, l_{n-i-1},
l_{n-i+2}, \dots, l_n).\end{eqnarray} Hence, using (\ref{adding})
and (\ref{fine}), we obtain
\begin{eqnarray}\label{fine1}
\qquad\quad  2 S_i(\ell) + M_i(\ell) \, \leq \, \left[2
S_{i+1}(\ell) + M_{i+1}(\ell)\right] + \left[2 S_{i-1}(\tilde
\ell) + M_{i-1}(\tilde \ell)\right].
\end{eqnarray}

 By our inductive hypothesis,
\begin{eqnarray*}
2\cdot S_{i+1}(\ell) + M_{i+1}(\ell) \, \leq\,  2^{n-i-1} - \sum_{j=r-i}^r \left(\begin{array}{c} n-i-1\\
j\end{array}\right) \, \\ \\
=\, 2^{n-i-1} - \sum_{j=r-i+1}^{r-1}
\left(\begin{array}{c} n-i-1 \\ j-1\end{array}\right) +\left(\begin{array}{c} n-i\\
r\end{array}\right)\end{eqnarray*} and
\begin{eqnarray*}
2 S_{i-1}(\tilde \ell) + M_{i-1}(\tilde \ell) \, \leq\,  2^{n-i-1} - \sum_{j=r-i+1}^{r-1} \left(\begin{array}{c} n-i-1\\
j\end{array}\right)
\end{eqnarray*}
Adding the last two inequalities and taking into account
(\ref{fine1}) we obtain
\begin{eqnarray*}
2 S_i(\ell) + M_i(\ell) \, &\leq&\,  2^{n-i} -\sum_{j=r-i+1}^{r-1}
\left[\left(
\begin{array}{c} n-i-1\\ j\end{array}
\right) + \left(
\begin{array}{c} n-i-1\\ j-1\end{array}
\right) \right] \\ \\ &+&\left(\begin{array}{c} n-i\\ r\end{array}  \right)=
2^{n-i} - \sum_{j=r-i+1}^r \left(
\begin{array}{c}
n-i\\ j
\end{array}
\right).
\end{eqnarray*}
This completes the proof of Theorem \ref{upper}.

To prove Theorem \ref{upper1} we assume that $n$ is even,
$n=2r+2$, and $\ell=(l_1, \dots, l_n)$ is a generic length vector
where $l_1\leq l_2\leq \dots\leq l_n$. We replace the inductive
hypothesis (\ref{intermid}) by

\begin{eqnarray}\label{intermid1}\quad\quad  2\cdot
S_i(\ell) \, \leq\,
2^{n-i}\, - \, 2\cdot \sum_{j=r-i+1}^r \left(\begin{array}{c} n-i-1\\
j\end{array}\right)\end{eqnarray} for $1\leq i\leq r+1.$ For $i=1$ inequality
(\ref{intermid1}) gives the desired inequality
$$2\cdot S_1(\ell)\, \leq\,  2^{n-1} - 2\cdot \left(\begin{array}{c} 2r\\
r\end{array}\right)\, =\, B'_n,$$ compare (\ref{sbmeven}). For
$i=r+1$ inequality (\ref{intermid1}) states $S_{r+1}(\ell) \leq 0$
which is obviously correct; this statement will be the base of
induction. To perform the step of induction we use inequalities
(\ref{adding}) and (\ref{fine}) which are valid in the case of $n$
even as well. We find
\begin{eqnarray*}
2\cdot S_{i+1}(\ell)\, \leq \, 2^{n-i-1} - 2\cdot \sum_{j=r-i}^r \left(
\begin{array}{c} n-i-2\\ j\end{array}
\right)
\end{eqnarray*}
and
\begin{eqnarray*}
2\cdot S'_i(\ell) \, \leq \, 2^{n-i-1} -2\cdot \sum_{j=r-i+1}^{r-1} \left(
\begin{array}{c}
n-i-2\\ j
\end{array}
\right)
\end{eqnarray*}
(both by the induction hypothesis) and adding the last two inequalities, using
(\ref{adding}), and performing transformations similar to the odd case, we
obtain (\ref{intermid1}).

This completes the proof of Theorem \ref{upper1}.

\end{document}